\documentclass[reqno,twoside,12pt]{amsart}

\hoffset - 1.5cm

\usepackage{amsmath,amssymb}
\usepackage{wasysym}
\usepackage{xcolor}

\newtheorem{definition}{Definition}[section]
\newtheorem{theorem}[definition]{Theorem}

\newtheorem{corollary}[definition]{Corollary}
\newtheorem{remark}[definition]{Remark}

\title{On the Green function to the Poisson and the Helmholtz equations on the $n$-dimensional unit sphere}
\author{Ilona Iglewska-Nowak}\address{West Pomeranian University of Technology, Department of Mathematics, al. Piastów 17, PL-70-310 Szczecin, Poland, ORCID 0000--0002--1938--8055 } \email{iiglewskanowak@zut.edu.pl}

\begin{document}

\begin{abstract}
A new method is presented to obtain a closed form of the generalized Green function to the Poisson and the Helmholtz equations on the $n$-dimensional unit sphere.
\end{abstract}

\keywords{spherical wavelets, $n$-spheres, PDE, Poisson equation, Helmholtz equation}
\subjclass[2010]{42C40, 42B37}

\date{\today}
\maketitle

\section{Introduction}

In the recently published paper~\cite{INS23}, written by Piotr Stefaniak and me, a formula for the Green function to the Poisson and the Helmholtz equations on the $n$-dimensional unit sphere was derived. The Green function is given as a series of the Gegenbauer polynomials. For the Poisson equation and some special values of the Helmholtz equation a closed form is given. In the present paper, a simpler method is derived to sum up the series. Consequently, computations are faster, and in the case of the Helmholtz equation a closed form of the Green function is obtained for a wider range of indices than in~\cite{INS23}.

\section{Preliminaries}\label{sec:preliminaries}

The present paper is strongly based on~\cite{INS23} and I refer to that publication for details on spherical nomenclature. For the purposes of the present investigation, the following is necessary.

\begin {enumerate}

\item The Poisson kernel for the $n$-dimensional sphere is given by
\begin{equation}\label{eq:Poisson_kernel}\begin{split}
p_r(y)&=\frac{1}{\Sigma_n}\sum_{l=0}^\infty r^l\,\frac{\lambda+l}{\lambda}\,\mathcal C_l^\lambda(\cos\vartheta)\\
&=\frac{1}{\Sigma_n}\frac{1-r^2}{(1-2r\cos\vartheta+r^2)^{(n+1)/2}},\quad r=e^{-\rho},\,y_1=\cos\vartheta,
\end{split}\end{equation}
where $\mathcal C_\kappa^K$ are the Gegenbauer polynomials of degree~$\kappa$ and order~$K$.

\item \cite[Theorem~1]{INS23} \begin{enumerate}\item Suppose that $f\in\mathcal C(\mathcal S^n)$ and $u\in\mathcal C^2(\mathcal S^n)$ satisfy
\begin{equation}\label{eq:Helmholtz}
\Delta^{\!\ast} u+au=f,
\end{equation}
where $a\in\mathbb R\setminus\{l(n+l-1),\,l\in\mathbb N_0\}$. Then,
\begin{equation}\label{eq:solution_u}
u=f\ast G,
\end{equation}
for
\begin{equation}\label{eq:Green_function_Helmholtz}
G=\sum_{l=0}^\infty\frac{1}{a-l(n+l-1)}\,\frac{\lambda+l}{\lambda}\,\mathcal C_l^\lambda.
\end{equation}
\item Suppose that $f\in\mathcal C(\mathcal S^n)$, $u\in\mathcal C^2(\mathcal S^n)$ are such that $\int f=\int u=0$ and satisfy
$$
\Delta^{\!\ast} u=f.
$$
Then,
\begin{equation}\label{eq:solution_u}
u=f\ast G,
\end{equation}
for
\begin{equation}\label{eq:Green_function_Poisson}
G=\sum_{l=1}^\infty\frac{-1}{l(n+l-1)}\,\frac{\lambda+l}{\lambda}\,\mathcal C_l^\lambda.
\end{equation}\end{enumerate}

\end{enumerate}

\newpage

\section{The main theorem}\label{sec:theorem}

In this section, a formula is given to sum up the series~\eqref{eq:Green_function_Helmholtz} or~\eqref{eq:Green_function_Poisson}. The advantage with respect to~\cite[Theorem~2]{INS23} is that the number of integrations is reduced to one.

\begin{theorem}\label{thm:Green_explicit} Let $n\in\mathbb N$, $n\geq2$, be fixed, $\lambda=\frac{n-1}{2}$, and suppose that $a=L(n+L-1)$ for some $L\in\mathbb R\setminus\mathbb Z$, $L\ne\frac{1-n}{2}$. Let  $L_0:=\max\{\left[L\right],\left[-n-L+1\right]\}$, where $\left[x\right]$ stays for the biggest integer less than or equal to~$x$. Denote by~$G$ the function
\begin{equation}\label{eq:Green}
G:=\sum_{l=0}^\infty\frac{1}{a-l(n+l-1)}\frac{\lambda+l}{\lambda}\,\mathcal C_l^\lambda.
\end{equation}
Then
\begin{align}
G&(\cos\vartheta)=\sum_{l=0}^{L_0}\frac{1}{a-l(n+l-1)}\cdot\frac{\lambda+l}{\lambda}\,\mathcal C_l^\lambda(\cos\vartheta)\label{eq:GL}\\
&+\frac{1}{n+2L-1}\int_0^1\left(r^{n+L-2}-r^{-L-1}\right)
   \left(\Sigma_n\,p_r(\cos\vartheta)-\sum_{l=0}^{L_0}r^l\,\frac{\lambda+l}{\lambda}\,\mathcal C_l^\lambda(\cos\vartheta)\right)\,dr.\notag
\end{align}
(If $L_0<0$, set $\sum_0^{L_0}=0$).
\end{theorem}

\begin{bfseries}Proof.\end{bfseries} Analogous to the proof of~\cite[Theorem~2]{INS23}. Note that
$$
\frac{1}{n+2L-1}\int_0^1\left(r^{n+L-2}-r^{-L-1}\right)\cdot r^l\,dr=\frac{1}{(L-l)(n+L+l-1)}
$$
for $l>L_0$.\hfill$\Box$

\begin{corollary}\label{cor:Green_explicit} Let $n\in\mathbb N$, $n\geq2$, be fixed, $\lambda=\frac{n-1}{2}$, and suppose that $a=L(n+L-1)$ for some $L\in\mathbb N_0$. Further, let~$G$ denote the function
\begin{equation}\label{eq:Green1}
G=\sum_{l=0,l\ne L}^\infty\frac{1}{a-l(n+l-1)}\frac{\lambda+l}{\lambda}\,\mathcal C_l^\lambda
\end{equation}
 Then
\begin{align}
G&(\cos\vartheta)=\sum_{l=0}^{L-1}\frac{1}{a-l(n+l-1)}\cdot\frac{\lambda+l}{\lambda}\,\mathcal C_l^\lambda(\cos\vartheta)\label{eq:G0}\\
&+\frac{1}{n+2L-1}\int_0^1\left(r^{n+L-2}-r^{-L-1}\right)
\left(\Sigma_n\,p_r(\cos\vartheta)-\sum_{l=0}^{L}r^l\,\frac{\lambda+l}{\lambda}\,\mathcal C_l^\lambda(\cos\vartheta)\right)\,dr.\notag
\end{align}
\end{corollary}

\begin{remark} If $a=L(n+L-1)$, then also $a=L^\prime(n+L^\prime-1)$ for $L^\prime=-n-L+1$. If $-\frac{(n-1)^2}{4}\leq a<0$, both~$L$ and~$L^\prime$ are negative (the case was not considered in \cite[Remark~3, p.~28]{INS23}). Choose $L=0$ for $a=0$ (Poisson equation) and~$L$ positive for positive~$a$ in order to apply Corollary~\ref{cor:Green_explicit}.
\end{remark}

\section{Closed forms of the Green functions}

\begin{theorem}
Let $a=L(n+L-1)$ for $L\in(1-n,\frac{1-n}{2})\cup(\frac{1-n}{2},0)$. Then
\begin{equation}\label{eq:G_sum_of_Appell}\begin{split}
G(\cos\vartheta)&=\frac{1}{n+2L-1}\cdot\left[\frac{F_1(n+L-1;\lambda+1,\lambda+1;n+L;e^{i\vartheta},e^{-i\vartheta})}{n+L-1}\right.\\
&-\frac{F_1(n+L+1;\lambda+1,\lambda+1;n+L+2;e^{i\vartheta},e^{-i\vartheta})}{n+L+1}\\
&-\frac{F_1(-L;\lambda+1,\lambda+1;-L+1;e^{i\vartheta},e^{-i\vartheta})}{-L}\\
&+\left.\frac{F_1(-L+2;\lambda+1,\lambda+1;-L+3;e^{i\vartheta},e^{-i\vartheta})}{-L+2}\right].
\end{split}\end{equation}
\end{theorem}

\begin{bfseries}Proof.\end{bfseries} Since $1-n<L<0$, $L_0$ defined in Theorem~\ref{thm:Green_explicit} is less than or equal to~$-1$ and
\begin{align}
G(\cos\vartheta)&=\frac{1}{n+2L-1}\int_0^1\left(r^{n+L-2}-r^{-L-1}\right)\frac{1-r^2}{(1-2r\cos\vartheta+r^2)^{\lambda+1}}\,dr\notag\\
&=\frac{1}{n+2L-1}\int_0^1\frac{r^{n+L-2}-r^{n+L}-r^{-L-1}+r^{-L+1}}{(1-2r\cos\vartheta+r^2)^{\lambda+1}}\,dr.\label{eq:G_big_negative_L}
\end{align}
All the exponents in the numerator of the integrand are greater than~$-1$.

According to~\cite[9.3(4)]{wB64},
\begin{equation}\label{eq:Appell_as_integral}\begin{split}
\int_0^1&u^{\alpha-1}(1-u)^{\gamma-\alpha-1}(1-ux)^{-\beta}(1-uy)^{-\beta'}du\\
&=\frac{\Gamma(\gamma-\alpha)\Gamma(\alpha)}{\Gamma(\gamma)}F_1(\alpha;\beta,\beta';\gamma;x,y),
\end{split}\end{equation}
where~$F_1$ is the Appell $F_1$-function and $0<\text{Re }\alpha<\text{Re }\gamma$. Since
$$
1-2r\cos\vartheta+r^2=(1-re^{i\vartheta})(1-re^{-i\vartheta}),
$$
we obtain from~\eqref{eq:Appell_as_integral}
$$
\int_0^1\frac{r^\alpha}{(1-2r\cos\vartheta+r^2)^{\lambda+1}}\,dr=\frac{F_1(\alpha;\lambda+1,\lambda+1;\alpha+1;e^{i\vartheta},e^{-i\vartheta})}{\alpha}.
$$
With this formula, \eqref{eq:G_sum_of_Appell} follows immediately from~\eqref{eq:G_big_negative_L}.\hfill$\Box$

\begin{remark}\begin{enumerate}
\item If~$n$ is even and~$L$ integer, the integrand in~\eqref{eq:G_big_negative_L} is a rational function of~$r$ and $\sqrt{1-2r\cos\vartheta+r^2}$. In this case, the function can be integrated with Euler substitution.
\item If $n$ is odd and $L$ rational, $L=-p/q$, the integrand in~\eqref{eq:G_big_negative_L} is a rational function of~$r^{1/q}$. Also in this case, another (than~\eqref{eq:Appell_as_integral}) integration rules are applicable and and one obtains another representation of~\eqref{eq:G_big_negative_L}.
\end{enumerate}\end{remark}

\end{document}